\documentclass[11pt]{amsart}
\usepackage{graphicx}
\usepackage{amsmath,amssymb,amsfonts,amscd,latexsym}

\newcommand{\labell}[1] {\label{#1}}

\newcommand{\1}{{{\mathchoice {\rm 1\mskip-4mu l} {\rm 1\mskip-4mu l}
{\rm 1\mskip-4.5mu l} {\rm 1\mskip-5mu l}}}}

\newcommand{\Cal}{{\rm Cal}}
\newcommand{\vol}{{\rm vol}}

\newcommand{\Tf}{{\Tilde f}}
\newcommand{\TQ}{{\Tilde Q}}

\newcommand\ft {{\mathfrak t}}
\newcommand{\Po}{{{\mathcal P}_\om}}

\newcommand{\Ver}{{\rm Vert}}

\newcommand{\less}{{\smallsetminus}}

\newcommand{\bla}{{\bigl\langle}}
\newcommand{\bra}{{\bigl\rangle}}

\newcommand{\THam}{{\Tilde\Ham}}

\newcommand{\p}{{\partial}}
\newcommand{\al}{{\alpha}}

\newcommand{\be}{{\beta}}

\newcommand{\Om}{{\Omega}}
\newcommand{\om}{{\omega}}
\newcommand{\eps}{{\varepsilon}}
\newcommand{\de}{{\delta}}
\newcommand{\De}{{\Delta}}
\newcommand{\ga}{{\gamma}}
\newcommand{\Ga}{{\Gamma}}
\newcommand{\io}{{\iota}}
\newcommand{\ka}{{\kappa}}
\newcommand{\la}{{\lambda}}
\newcommand{\La}{{\Lambda}}
\newcommand{\si}{{\sigma}}
\newcommand{\Si}{{\Sigma}}

\newcommand{\Isom}{{\rm Isom}}

\newcommand{\Mm}{{\mathcal M}}
\newcommand{\Ss}{{\mathcal S}}
\newcommand{\Ll}{{\mathcal L}}

\newcommand{\ov}{\overline}

\newcommand{\Tpsi}{{\Tilde \psi}}
\renewcommand{\Tilde}{\widetilde}
\newcommand{\TM}{{\Tilde M}}
\newcommand{\Tphi}{{\Tilde \phi}}

\newcommand{\TP}{{\Tilde P}}

\newcommand{\Tom}{{\Tilde \om}}
\newcommand{\TOm}{{\Tilde \Om}}

\newcommand{\PP}{{\mathbb P}}
\newcommand{\Q}{{\mathbb Q}}
\newcommand{\R}{{\mathbb R}}

\newcommand{\C}{{\mathbb C}}
\newcommand{\SO}{{\rm SO}}

\newcommand{\Z}{{\mathbb Z}}

\newcommand{\Ham}{{\rm Ham}}
\newcommand{\Symp}{{\rm Symp}}

\newcommand{\Flux}{{\rm Flux}}
\newcommand{\Diff}{{\rm Diff}}
\newcommand{\Vol}{{\rm Vol}}

\newcommand{\Aa}{{\mathcal A}}
\newcommand{\Nn}{{\mathcal N}}
\newcommand{\Pp}{{\mathcal P}}

\newcommand{\ev}{{\rm ev}}
\newcommand{\SSS}{{\smallskip}}
\newcommand{\QED}{{\hfill $\Box$\MS}}

\newtheorem{theorem}{Theorem}[section]

\newtheorem{cor}[theorem]{Corollary}
\newtheorem{lemma}[theorem]{Lemma}

\newtheorem{prop}[theorem]{Proposition}

\newtheorem{rmk}[theorem]{Remark}
\newtheorem{quest}[theorem]{Question}
\numberwithin{figure}{section}
\numberwithin{equation}{section}
\numberwithin{table}{section}

\newcommand{\MS}{{\medskip}}

\newcommand{\NI}{{\noindent}}

\begin{document}

 \title{Loops in the Hamiltonian group: a survey}
 \author{Dusa McDuff}\thanks{partially supported by NSF grant DMS 0604769. Some of this material formed the basis of a talk given at
 the AMS summer conference on Hamiltonian dynamics in Snowbird, UT, July 2007, which was also supported by the NSF}
\address{Department of Mathematics,
Barnard College, Columbia University, New York, NY 10027-6598, USA.}
\email{dusa@math.columbia.edu}
\date{17 November 2007, revised 18 January 2009}

\begin{abstract}
This note describes some recent results about the
homotopy properties of Hamiltonian loops 
in various manifolds, including toric manifolds and 
one point blow ups.
We describe conditions
 under which a circle action does not contract in the Hamiltonian group, and construct an example of a 
  loop $\ga$  of 
 diffeomorphisms of a symplectic manifold $M$   with the property that 
 none of the loops smoothly isotopic to $\ga$ preserve any symplectic form on $M$. We also discuss some new conditions under which
  the Hamiltonian group has infinite Hofer diameter.
Some of the methods used are classical (Weinstein's action homomorphism and volume calculations), while others use quantum methods (the Seidel representation and spectral invariants).
\end{abstract}

\maketitle





\section{Introduction}

Let  $(M,\om)$ be a symplectic manifold that is closed, i.e. compact and without boundary. We denote its group of diffeomorphisms by $\Diff: = \Diff\, M$ and by $\Symp: = \Symp(M,\om)$ 
the subgroup  of symplectomorphisms, i.e. diffeomorphisms
 that preserve the symplectic form.   Its identity component $\Symp_0: = \Symp_0(M,\om)$ has an important normal subgroup
 $\Ham: = \Ham(M,\om)$ consisting of all symplectomorphisms $\phi\in \Symp_0$ with zero flux, or equivalently, of all time $1$ maps $\phi_1^H$  of Hamiltonian flows $\phi_t^H$, where $H: M\times S^1\to \R$ is a (smooth) time dependent  function. Basic information about 
 these groups may be found in
 \cite{MS1,MS2,Pbk} and the survey articles \cite{Mcgr,Mcox}.
 Note that $\Ham=\Symp_0 $ when $H^1(M;\R) = 0$.
 
This survey is mostly concerned with questions about based Hamiltonian loops, i.e. smooth paths  $\{\phi_t\}_{0\le t\le 1}$ in the Hamiltonian group $\Ham$ for which $\phi_1=\phi_0=id$. Each such loop is 
the flow of some time dependent Hamiltonian $H_t,$ and may be reparametrized so that $H_0 = H_1$ and the induced map 
$H: M\times S^1\to \R$ is smooth.   If this generating Hamiltonian is time independent then its flow $\phi^H_t, t\in S^1$, is a subgroup
of $\Ham$ isomorphic to $S^1$. (The function $H:M\to \R$ is then called the {\it moment map}.)   Although these loops are the easiest to understand,
 there are still many unsolved questions about them.

Our first group of questions concerns the structure of
$\pi_1(\Symp)$ and
$\pi_1(\Ham)$.   Note that when $\dim \, M = 2$, Moser's homotopy argument implies that the symplectomorphism group of $M$ is homotopy equivalent to its group of orientation preserving
diffeomorphisms.  Thus the homotopy type of the groups $\Symp$ and $\Ham$ are known.  $\Ham$ is contractible unless $M=S^2$, in which case both it and $\Symp_0$ are 
 homotopic to  $\SO(3)$; while in higher genus $\Symp_0$ is contractible except in the case of the torus, when it is homotopic to the torus $\R^2/
 \Z^2$.
The group $\pi_0(\Symp)$ is the well-known mapping class group.
 Thus the questions listed below are not interesting in this case.  
 
 The case when
 $\dim M=4$ is also moderately well understood.  In particular, every $4$-manifold with a Hamiltonian $S^1$ action is the blow up of a
 rational or ruled $4$-manifold. (This was first proved by Audin \cite{Aud} and Ahara--Hattori \cite{AH}; see also Karshon \cite{Kar}.)  Moreover the homotopy type of 
 $\Ham$ is understood when $M = \C P^2$ or $S^2\times S^2$ or a one point blow up of such; see for example Gromov \cite{G}, 
 Abreu--McDuff \cite{AM}, Abreu--Granja--Kitchloo \cite{AGK} and Lalonde--Pinsonnault \cite{LALPi}.
  For work on nonHamiltonian $S^1$ actions in dimensions $4$ and above see Bouyakoub~\cite{Bou}, Duistermaat--Pelayo~\cite{DP}
   and Pelayo~\cite{Pel}.    There is information on  smooth circle actions in $4$-dimensions, see Fintushel~\cite{Fin} 
   and Baldridge~\cite{Bald}. 
But almost nothing is known about the diffeomorphism group of a manifold of dimension $\ge 4$; for example  it is not known
 whether 
$\Diff\, \C P^2$ is homotopy equivalent to the projective unitary group $PU(3)$ as is $\Symp(\C P^2)$. 
On the other hand it follows from \cite{AM} that
$\Ham(S^2\times S^2,\om)$ is not homotopy equivalent to 
$\Diff(S^2\times S^2)$ for any symplectic form $[\om]$.

\begin{quest}\labell{q1}  When does a circle subgroup $\ga$  of $\Symp$ represent a nonzero element in $\pi_1(\Symp)$,  or even one
 of infinite order?
\end{quest}

(Entov--Polterovich~\cite{EP} call circles of infinite order
{\it incompressible}.)

The first problem here  is to decide when
a loop is Hamiltonian, i.e. is in the kernel of the 
Flux homomorphism.  Recall from \cite{MS1} that Flux is defined on the universal cover $\Tilde\Symp_0$ of $\Symp_0$
by
\begin{equation}\labell{eq:flux}
\Flux: \Tilde\Symp_0\to H^1(M;\R),\quad
 \Tphi \mapsto \int_0^1 [\om(\dot\phi_t,\cdot)]\,dt,
\end{equation}
where $\Tphi = (\phi_1, \{\phi_t\}) \in \Tilde\Symp_0$.
The (symplectic) Flux group $\Ga_\om$ is defined to be the image of $\pi_1(\Symp) $ under Flux, so that there is an induced homomorphism
$$
\Flux:\pi_1(\Symp)\to H^1(M;\R)/\Ga_\om
$$ 
with kernel  $\Ham$.
Ono~\cite{Ono} recently proved that  $\Ga_\om$ is discrete.

 Unfortunately there seem to be no good techniques for 
 understanding when $\Flux(\ga)$ is trivial.   Since Hamiltonian $S^1$ actions always have fixed points at the critical points of the  moment map $H:M\to \R$, a first guess might be that every symplectic action with fixed  points is Hamiltonian.  However
McDuff~\cite{Mc88} shows that this is not the case except in dimension $4$. 
In fact, the following basic problem is still unsolved in dimensions $>4$.  

\begin{quest}\labell{q4}  Suppose that $S^1$ acts symplectically on the closed  symplectic manifold $(M,\om)$ with a finite but nonzero number of fixed points.  Is the action Hamiltonian?
\end{quest}

If the action is semifree (i.e. the stabilizer of a point is either the identity or the whole group) Tolman--Weitsman~\cite{TW} show that the answer to Question~\ref{q4} is affirmative by computing various equivariant cohomology classes.
Some other information on this question has been obtained by
  Feldman~\cite{Fel} and Pelayo--Tolman
  ~\cite{PelT}.  One might hope to extend 
  the Tolman--Weitsman result to semifree actions with more general conditions on the fixed point components, for example that they are simply connected; note that these cannot be arbitrary because of the example in \cite{Mc88} of a semifree but
   nonHamiltonian action on a $6$-manifold
   with fixed point sets that are $2$-tori. 
  
  In the current discussion we will largely ignore this problem, for the most part considering only Hamiltonian loops and their images in $\pi_1(\Ham)$.

\begin{quest} \labell{q2} To what extent are $\pi_1(\Ham)$ 
and $\pi_1(\Symp)$ generated by symplectic $S^1$ actions?
\end{quest}

This question is a measure of our ignorance.  $S^1$ actions do generate
 $\pi_1(\Symp)$ in very special cases
  such as $\C P^2$ or $\bigl(S^2\times S^2, pr_1^*(\si) + pr_2^*(\si)\bigr)$ 
  (note that the factors have equal area). Indeed
in these cases  $\Symp$ itself is known to have the homotopy type of a compact Lie group (see \cite{Mcgr}).
However, as we see below, this does not hold in general.

\begin{quest} \labell{q3} What can one say about the relation between 
$\pi_1(\Ham)$, $\pi_1(\Symp)$ and $\pi_1(\Diff)$?  For example, under what circumstances is the map
$\pi_1(\Symp)\to \pi_1(\Diff)$ injective or surjective?
\end{quest}

The symplectic Flux group $\Ga_\om$  is the quotient of $\pi_1(\Symp)$ by
$\pi_1(\Ham)$ and hence precisely measures their difference.
By K\c edra--Kotschick--Morita \cite{KKM}, this group vanishes in many cases.  Much of their paper in fact applies to the volume\footnote
{
The volume flux is defined by equation~(\ref{eq:flux}), but $\om$ 
should be understood as a volume form and the homomorphism takes 
values  in $H^{m-1}(M)$, where $m: = \dim M$. 
Accordingly, $\Ga_{\vol}$ is the
image of $\pi_1(\Diff_{\vol})$ under the flux.}  
flux group $\Ga_{\vol}$, which is in principle of a  more topological nature than $\Ga_\om$; it would be 
interesting to  find conditions
for the vanishing of $\Ga_\om$  that involve 
symplectic geometry at a deeper level.

In this note, we begin by describing 
 some classical methods for exploring 
 the above questions, the first based on properties of the {\it action functional $\Aa_H$} and the second using volume.  These methods give rather good  information in the following  cases:
\MS

\NI$\bullet$  Question~\ref{q1} for toric manifolds (see Corollary~\ref{cor:tor} below)

\NI$\bullet$   Questions~\ref{q2} and
~\ref{q3}  for pointwise blow ups $\TM$ of arbitrary symplectic manifolds 
$M$  (see Proposition~\ref{prop:csymp} and its corollaries).\MS

If $M$ is noncompact and one considers the group $\Ham^c: = \Ham^c\,M$ of compactly supported Hamiltonian symplectomorphisms of $M$,
then there is another classical homomorphism called the {\it Calabi homomorphism}:
$$
\Cal:\pi_1(\Ham^c)\to \R,\quad \ga\mapsto \frac1{n!}\int_0^1\Bigl(\int_M H_t\om^n\bigr) dt,
$$
where $H_t$ is the generating Hamiltonian for $\ga$, normalized to have compact support.  We explain briefly in Lemma~\ref{le:cal} why\SSS

\NI$\bullet$ the Calabi homomorphism need not vanish.\SSS

\NI
As we shall see, this question, though classical in origin,  is very closely related to questions about the Seidel representation in quantum homology.\MS

One might wonder if it is possible to get better information about the above questions by using  more modern (i.e. 
quantum) techniques.     In fact, Question~\ref{q1} first 
arose in McDuff--Slimowitz \cite{MSlim}, a paper that uses Floer theoretic techniques to study paths in $\Ham$ that are geodesic with respect to the Hofer norm.\footnote
{
This is defined in \S\ref{ss:Hof} below.} 
 An unexpected consequence of the ideas developed there is that semifree
Hamiltonian circle actions  do not contract in $\Ham$, though they might
have finite order (for example, a rotation of $S^2$.) 
 The main tool that has proved useful in this context is an extension of the action homomorphism due to Seidel~\cite{Sei}, that is called the {\it Seidel representation}; see \S\ref{ss:seid}.   This homomorphism assigns to every $\ga\in \pi_1(\Ham)$ a
unit (i.e. invertible element) $\Ss(\ga)$ in the (small) quantum homology $QH_*(M)$.
Corollary~\ref{cor:seid} gives some more results on the above questions obtained using $\Ss$.  Because $\Ss$ is usually very hard to calculate,
 the classical methods often work better in specific examples.  Nevertheless,
 $\Ss$ is a key tool in other contexts.
 
 One very interesting question is the following. Note that in two dimensions, the only symplectic manifold with a Hamiltonian $S^1$ action is $S^2$, 
 while $T^2$ has nonHamiltonian actions and higher genus surfaces $\Si_g$ have none.

\begin{quest}\labell{q6}  Is there a meaningful extension 
 of the classification of Riemann surfaces into
 spheres, tori  and higher genus  to higher dimensional symplectic manifolds?
 If so, is any aspect of it reflected in the properties of $\pi_1(\Ham)$?
 \end{quest}
 
This would be an analog of {\it minimal model} theory
 in algebraic geometry.  A first step, accomplished by 
 Ruan and his coworkers Hu and T.-J. Li \cite{R,HLR}, 
 is to understand what it means for two
  symplectic manifolds to be birationally equivalent. 
  Their results imply that a reasonable  class of manifolds to take as the analog of spheres are the {\it symplectically uniruled} manifolds.  These are the manifolds for which there is a nonzero genus zero Gromov--Witten invariant
$\bla a_1,a_2,\dots,a_k\bra^M_\be$ (for some $k\ge 1,$ $a_i\in H_*(M)$ and $\be\in H_2(M)$) with one of the constraints $a_i$ equal to a point. 
This class includes all projective manifolds that are uniruled in the sense of algebraic geometry.
In this case  $\om(\be)\ne0$ and $c_1(\be)\ne 0$.  At the other extreme
 are the {\it symplectically aspherical}  manifolds for which the restriction $\om|_{\pi_2(M)}$ of $[\om]$ to $\pi_2(M)$ vanishes, and possibly also (depending on the author) the restriction
$c_1|_{\pi_2(M)}$ of the first Chern class $c_1$ of $(M,\om)$.  These manifolds have no $J$-holomorphic curves at all, and hence all nontrivial (i.e. $\be
\ne 0$) Gromov--Witten invariants vanish.

To a first approximation, one can characterize symplectically 
uniruled manifolds in terms of 
their quantum homology $QH_*(M)$ in the following way.  

\begin{quote}
{\it If $(M,\om)$ is not uniruled then all invertible elements in $QH_{2n}(M)$ have the form $\1\otimes \la + x$ where $\la$ is invertible in the coefficient field $\La$ and $x\in H_{<2n}(M)\otimes \La$.}
\end{quote}

\NI
(This is nearly an iff statement and can be
improved to such: cf.  the appendix to \cite{Mcu}.
Here $\1$ denotes the fundamental class $[M]\in QH_*(M)$; it is the identity element in $QH_*(M)$.  The notation is explained in more detail below.)   Therefore another version of the previous question is:

\begin{quest}\labell{q7}  To what extent is the structure of the quantum homology ring $QH_*(M)$ reflected in the algebraic/topological/geometric structure of $\Ham$?
\end{quest}

For example, according to Polterovich~\cite{P2},

\begin{quote}{\it if $\om|_{\pi_2(M)} = 0$, $\Ham$ contains no elements of finite order except for $id$.}
\end{quote}

Does this property continue to hold if the condition on $\om$ is weakened
to the vanishing of all Gromov--Witten invariants?
Since $S^1$-manifolds certainly support nontrivial symplectomorphisms of finite order (such as a half turn), it is natural to ask whether these manifolds have nontrivial
Gromov--Witten invariants.  
Since the Seidel element $\Ss(\ga)$ is a unit in quantum homology
one might also expect to see traces of the uniruled/non-uniruled  dichotomy in its properties. This was the guiding idea in my recent proof
\cite{Mcu} that\MS

\begin{quote} {\it every closed symplectic manifold that supports a Hamiltonian $S^1$ action is uniruled}.
\end{quote}

\NI
The argument in \cite{Mcu} applies
more generally  to manifolds with Hamiltonian loops that are nondegenerate 
 Hofer geodesics.  This opens up many interesting questions of a more dynamical flavor.

 Here we shall discuss the following basic (but unrelated) problem, which  is still open in many cases. 

\begin{quest}\labell{q8} Does 
$\Ham$ have infinite diameter with respect to the Hofer metric?
\end{quest}

\NI   
One expects the answer to be positive always.  However the proofs  for spheres (due to Polterovich~\cite{P}) and other Riemann surfaces (due to Lalonde--McDuff~\cite{LM0}) are very different. In fact, as noted by
Ostrover~\cite{Os}  
one can use the spectral invariants of
Schwarz~\cite{Sch} and Oh \cite{Oh,Oh1} (see also Usher~\cite{U}) to show that  the universal cover
$\THam$  of $\Ham$  always has infinite diameter with respect to the 
induced (pseudo)metric.  Therefore the question becomes: when does this result transfer down to $\Ham$?  Schwarz shows in~\cite{Sch} that this happens when
$\om$ and $c_1$ both vanish on $\pi_2(M)$.   I recently
extended his work in  \cite{Mcmdr}, showing that the asymptotic spectral invariants descend to $\Ham$ if, for example, all nontrivial genus zero
Gromov--Witten invariants vanish and ${\rm rank\,} H_2(M;\R)>1$.
 As we explain in \S\ref{ss:Hof} below, 
 Schwarz's argument 
 hinges on the properties of 
the Seidel elements $\Ss(\ga)$ of  $\ga\in \pi_1(\Ham)$.
Here we shall sketch a  different extension of his result.
In particular we show:

\begin{quote}
{\it $\Ham$ has infinite diameter when $M$ is a \lq\lq small" one point blow up of $\C P^2$.}
\end{quote}

This manifold $M$ is of course uniruled (and the spectral invariants do not descend).  To my knowledge, it is not yet known
whether $\Ham$ has infinite diameter for all one point 
blow ups of $\C P^2$, though it does for $\C P^2$ itself (and indeed for any $\C P^n$) by the results of Entov--Polterovich~\cite{EP1}.
For further results on this problem 
see McDuff \cite{Mcmdr}.

\section{Classical methods}\labell{s:class}

We consider two classical methods to detect elements in $\pi_1(\Ham)$, 
the first based on the action 
homomorphism and the second on considerations of volume.
One can also use homological methods as in K\c edra--McDuff~\cite{KM} but these give more information on the higher homotopy groups $\pi_k(\Ham), k>1$.
  
\subsection{The action homomorphism}\labell{ss:act}

First of all, we consider Weinstein's {\it action homomorphism}
$$
\Aa_\om:\pi_1(\Ham) \to \R/\Pp_\om,
$$
where $\Pp_\om: = \{\int_c\om|c\in \pi_2(M)\}\subset \R$ is the group of spherical periods of $[\om]$.  In defining $\Aa_\om$, we shall use the following sign conventions.  The flow $\{\phi^H_t\}_{t\in [0,1]}$ of a Hamiltonian $H:M\times S^1\to \R$ satisfies the equation
$$
\om(\dot{\phi}^H_t(p),\cdot) = -dH_t(\phi^H_t(p)),\quad p\in M,
$$
and the corresponding action functional on the space $\Ll_0M$ of contractible loops in $M$ is
$$
\Aa_H: {\Ll}_0M\to \R/\Pp_\om,\qquad \Aa_H(x): = -\int_{D^2}u^*\om
+\int_{S^1}H_t(x(t)) \,dt,
$$
where $u:D^2\to M$ is any extension of the (contractible) loop
$x:S^1\to M$, and $H_t$ is assumed to be {\it mean normalized}, i.e. $\int_MH_t\om^n = 0$ for all $t$. Note that any Hamiltonian $H_t$ can be normalized (without changing its flow) by subtracting a suitable {\it normalization constant} $c_t: = {\int_MH_t\om^n}\,/{\int_M\om^n}$.

One very important property of $\Aa_H$ is
that its critical points are the closed orbits $x(t): = \phi_t^H(p)$ of the corresponding flow.  This is classical: the orbits of a Hamiltonian flow minimize (or, more correctly, are critical points of) the action. 
If $\{\phi_t^H\}_{t\in S^1}$ is a loop,
 then all its orbits  are closed. Moreover they are contractible.\footnote{
This folk theorem is proved, for example,
 in \cite[Ch~9.1]{MS2} in the case when $(M,\om)$ is semipositive.  The general case follows from the proof of the Arnold conjecture: see \cite{FO,LT}. A simpler proof is sketched in \cite{Mcox}.
Since all known proofs use quantum methods the existence of the homomorphism $\Aa_\om$ is not entirely \lq\lq classical" in the case when $H^1(M)\ne \{0\} $.  However, toric manifolds are simply connected, and so our discussion of that case is \lq\lq classical".}
Hence these orbits all have the same action.  It is not hard to see that this value depends only on the homotopy class $\ga: = [\{\phi_t^H\}]$ of the loop. Call it $\Aa_\om(\ga)$.  We therefore get a map $\Aa_\om:\pi_1(\Ham) \to \R/\Pp_\om$, which is easily seen to be a homomorphism.

If the loop $\{\phi_t\}_{t\in S^1}$ is  the circle subgroup
$\ga_K$ of $\Ham(M,\om)$ generated by the  mean 
normalized Hamiltonian $K \colon M\to \R$ then 
we may take $p$ to be a critical point of $K$ and $u$ to be the constant disc.  Hence 
$\Aa_{\om}(\ga_K)$ is the image in $\R/\Po$ of 
any critical value of $K$.
This is well defined because the difference $K(p)-K(q)$ 
between two critical values is the integral of $\om$ over the 
$2$-sphere formed by rotating an arc from $p$ to $q$ 
by the $S^1$-action.

Polterovich noticed that when $(M,\om)$ is spherically monotone 
(i.e. there is a constant $\ka$ such that $\ka [\om]$ and $c_1(M)$ induce the same homomorphism $\pi_2(M)\to \R$) then one can use 
the Maslov index of a fixed point to get rid of the 
indeterminacy of $\Aa_\om$, hence lifting it to 
 the {\it combined Action--Maslov homomorphism} $\pi_1(\Ham)\to\R$; cf.   
\cite{Parn,EP}.

In general, the action homomorphism is hard to calculate because 
there are few good ways of understanding
the normalization constant for an arbitrary Hamiltonian.  However,
as pointed out in McDuff--Tolman~\cite{MT2}, this calculation is possible in the toric case and, as we now explain, one gets quite useful information from it. 

\MS

\NI {\bf The toric case.}

A $2n$-dimensional symplectic manifold is said to be
{\it toric} if it supports an effective Hamiltonian action of the standard
torus $T^n$.  
The image of the corresponding moment map $\Phi:M\to \R^n$ is a convex polytope $\De$ (called the {\it moment polytope})
 that, as shown by Delzant, completely determines the initial symplectic manifold which we therefore denote $(M_\De,\om_\De)$.    
These manifolds have a  compatible complex structure $J_\De$, also determined by $\De$, and hence are K\"ahler; cf. Guillemin \cite{Gui}.\footnote
{
In fact there are many  compatible complex structures on $M_\De$.  By $J_\De$ we mean the one obtained by thinking of $M_\De$ as a symplectic quotient of $\C^N$, where $N$ is the number of (nonempty) facets of $\De$; see for example \cite[Ch~11.3]{MS2}.
}
  Denote by $\Isom_0\,
M_\De$ the identity component of the corresponding Lie group of K\"ahler isometries. (By \cite{MT2}, this is a maximal 
compact connected subgroup of $\Symp\,M_\De$.) 
Then $T^n\subseteq \Isom_0\,M_\De$, and it is natural to try to 
 understand the resulting homomorphisms
\begin{equation}\labell{eq:iso}
\pi_1(T^n)\to \pi_1(\Isom_0\,M_\De) \to \pi_1(\Ham\,M_\De).
\end{equation}
For example, if $M = \C P^n$ (complex projective space), $\De$ is an $n$-simplex and $\Isom_0\,M_\De$ is  the projective unitary group $= PU(n+1)$. 
Hence the first map has image $\Z/(n+1)\Z$.    Seidel pointed out in his thesis \cite{Sei} that this finite subgroup $\Z/(n+1)\Z$ injects into $\pi_1(\Ham\,\C P^n)$.  We will prove this here by calculating the action homomorphism. \footnote{
Although our proof looks different from Seidel's, it is essentially the same; cf. the description (given below) of $\Aa_\om(\ga)$ in terms of the fibration $P_\ga\to S^2$.}
 
Denote by $\ft$ the Lie algebra of $T^n$, and by $\ell$ its
integral lattice.  Each vector $H\in \ell$ exponentiates to 
a circle subgroup $\ga_H\subset T^n$.  Strictly speaking, the moment map should be considered to take values in the dual $\ft^*$ of the Lie algebra
$\ft$, and its definition 
implies that  the corresponding  $S^1$ action on $M$ is generated by the function
 $$
p\mapsto \bla H, \Phi(p)\bra, \qquad p\in M,
 $$
 where $\bla\cdot,\cdot\bra$ denotes the natural pairing between $\ft$ and its dual $\ft^*\supset {\rm Im}\,\Phi$.
 The corresponding mean normalized Hamiltonian is 
$$
K: =  \langle H,  \Phi  \rangle  - \langle H, c_\Delta \rangle,
$$ 
where $c_\Delta$ denotes the center of mass of $\Delta.$
Since vertices of $\Delta$ correspond to fixed points in $M_\Delta$,
this implies that
$$
\Aa_{\om}(\ga_H) \;=\;  \langle H ,   v  \rangle - \langle H, c_{\Delta} \rangle 
\;\;\in \;\;
\R/\Po,
$$
where $v$ is any vertex of $\Delta$.

To go further, recall that the moment polytope $\De$ is {\it rational}, that is,
the outward conormals  to the facets $F_i, 1\le i\le N$, of $\De$ are 
rational and so have unique primitive representatives 
$\eta_i\in \ell\subset\ft, 1\le i\le N.$  (This implies that slopes of the edges of $\De$ are rational, but its vertices need not be.) Thus $\De$ may be described as the solution set of a system of linear  inequalities:
$$
\De= \bigcap_{i=1}^N\{x\in \ft^*: \bla \eta_i,x\bra\le \ka_i\}.
$$
The constants $ (\ka_1,\dots,\ka_N)$ are called the {\it support numbers of $\De$} and determine the symplectic form $\om_\De$. 
They may be slightly varied without changing the diffeomorphism type of $M_\De$. We shall write $\ka: = (\ka_1,\dots,\ka_N)$.
Thus we may think of  
the moment polytope $\De: = \De(\ka)$ and  its center of gravity $c_\De(\ka)$ as functions of  $\ka$.
An element $H\in \ell$  is said to be a {\it mass linear function} on $\De$
if  the quantity $\bla H,c_\De(\ka)\bra$ depends  linearly on $\ka$. 

Here is a foundational result from \cite{MT2}.  The proof that the $\al_i$ are integers uses the fact that $M_\De$ is smooth, i.e. for each vertex of $\De$ the conormals of the facets meeting at that point form a basis for the integer lattice $\ell$.

\begin{prop}\labell{prop:symp} 
Let $H \in \ell$.  If $\ga_H$   
contracts in $\pi_1\bigl(\Ham(M_\Delta, \om_{\De(\ka_0)})\bigr)$ then $H$ is mass linear. More precisely, there are integers $\al_i$ such that
for all $\ka$ sufficiently near $\ka_0$
$$
\langle H, c_\Delta(\kappa) \rangle  =  \sum \al_i \kappa_i.
$$
\end{prop}

\begin{proof}[Sketch of proof]
Clearly if $\ga_H$ vanishes in 
$\pi_1(\Ham\,M_\Delta)$, then $\Aa_{\om}(\ga_H) = 0$.
It is immediate from Moser's homotopy argument that  
if $\ga_H$ contracts in $\Ham(M_\Delta,\om_\Delta)$, 
then it also contracts for all sufficiently small perturbations of 
$\omega_\Delta$.   
Since varying $\kappa$ corresponds to varying the
symplectic form  on $M_\Delta$, for any vertex $v$ of $\De$
the image of $\ga_H$ under the action homomorphism
$$
\Aa_{\om(\ka)}(\ga_H) 
= \langle H, v \rangle  - \langle H,c_{\Delta(\ka)} \rangle
$$
must lie in ${\mathcal P}_{\om(\ka)}$ for all $\ka$ 
sufficiently close to $\ka_0\in \R^N$.
But ${\mathcal P}_{\om(\ka)}$ is generated by the (affine) lengths\footnote
{
Since $\De$ is rational its edges $e$ have rational slopes and so for
 each $e$ there is
an integral affine transformation of $(\ft^*,\ell^*) \cong (\R^{n},\Z^n)$ taking it to the $x_1$ axis.  
The affine length of $e$ is defined to be the usual length of its image on the $x_1$-axis.}
 of the  edges
of the polytope $\Delta$, and so
it is a finitely generated subgroup of $\R$ 
whose generators are linear functions of the $\kappa_i$ with
integer coefficients.
Similarly, for each vertex $v=v(\ka)$ of $\De(\ka)$ 
the function $ \langle H, v \rangle$ is  linear with respect to 
 the $\kappa_i$ with integer coefficients.
Since the function $\ka\mapsto \Aa_{\om(\ka)}(\ga_H)$ is 
continuous, it follows that the 
function  $\ka \mapsto \langle H,c_{\Delta(\ka)} \rangle$
is also a linear function of the $\kappa_i$ with integer coefficients
 as $\ka$ varies in some open set.  
\end{proof}

We now can sketch the proof of the result mentioned just 
after equation (\ref{eq:iso}) that $\pi_1(PU(n+1))$ injects into 
$\pi_1(\Ham\, \C P^n)$.

\begin{prop}  If $M_\De = \C P^n$, the map $\pi_1(\Isom_0\,M_\De)\to \pi_1(\Ham\,M_\De)$ is injective.
\end{prop}
\begin{proof}   Since $\De$ is a simplex, we can
take it to be $\De: = \{x\in \R^n|x_i\ge 0,\sum x_i = 1\}$ with
edges of affine length $1$. Hence $\Po=\Z$ and the 
 center of gravity is $(\frac 1{n+1},\dots,\frac 1{n+1})$.  Next, note that $\pi_1(\Isom_0\, M_\De)$ is generated by the circle action
$$
\ga: [z_0:z_1:\dots:z_n]\mapsto [e^{-2\pi i\theta} z_0:z_1:\dots: z_n].
$$
The corresponding vector $H\in \ft$ is 
 one of the outward conormals of the moment polytope $\De$.  One can easily check that 
$\langle H,c_{\Delta} \rangle = -\frac 1{n+1}\in \R/\Po=\R/\Z$, which has order $n+1$.  
Hence the order of $\ga$ in $\pi_1(\Ham)$ is divisible by 
$n+1$, and so the map must be injective.\end{proof}

It turns out that most moment polytopes
 do not admit mass linear $H$.
 More precisely, we call a facet {\it pervasive} if it  meets all the other facets and {\it flat} if the conormals of all the facets meeting 
 it lie in a hyperplane in $\ft$. Note that one can 
 destroy pervasive and flat facets by suitable blow ups.

\begin{prop}[\cite{MT2}]\labell{prop:tor}  Suppose that $\De$ is a moment polytope with no pervasive or flat facets.  Then
every   mass linear function  $H$ on $\De$ vanishes.
\end{prop}

\begin{cor} \labell{cor:tor}  If $\De$ has no pervasive or flat facets,
 the map $\pi_1(T^n)\to 
\pi_1(\Ham\,M_\De)$ is injective. 
\end{cor}

On the other hand we also show:

\begin{prop} [\cite{MT2}]\labell{prop:tor2} The only toric manifolds for which  the image of $\pi_1(T^n)$  in $\pi_1(\Ham\, M_\De)$ is finite are products of projective spaces.
\end{prop}

 It is also easy to characterize mass linear $H$ corresponding to circles $\ga_H$ that contract 
in $\Isom_0\,M_\De$. (Such $H$ are called {\it inessential}.  The other mass linear functions are called {\it essential}.)  The papers~\cite{MT2,MT3} classify all pairs $(\De, H)$ consisting of  a moment 
polytope $\De$ of dimension less than or equal to four together with an essential
 mass linear  $H$ on it.  It turns out that there are many 
 interesting families of examples 
when $n=4$.  However when $n\le 3$ there is only one.

\begin{prop}  Let $(M_\De, \om_\De)$ be a  toric manifold 
of dimension $2n\le 6$ such that $\pi_1(\Isom_0\,M_\De)\to \pi_1(\Ham\,M_\De)$ is 
not injective.  Then $M$ is a $\C P^2$-bundle over $\C P^1$.
\end{prop}

\subsection{Applications using Volume}\labell{ss:vol}

This section is based on the paper \cite{Mcbl}, which uses various methods, both classical and quantum, to explore the homotopy groups of $\Ham\,\TM$, where $\TM$ is the one point blow up of $M$.
We shall concentrate here on questions concerning $\pi_1(\Ham\,\TM)$.
Our arguments  do not use the geometry provided by the form $\om$, but just the fact that its cohomology class $a: = [\om]$ has $a^n\ne 0$.
Thus in this section we will suppose that $(M,a)$ is a {\it $c$-symplectic manifold} (short for cohomologically symplectic), which simply means that $a\in H^2(M)$ has the property that $a^n\ne 0$, where, as usual, $2n: = \dim M$.   Note that such a manifold is oriented, and so for each point $p\in M$ has a well defined blow up $\TM$ at $p$ obtained by choosing a complex structure near $p$ that is compatible with the orientation and then performing the usual complex blow up at $p$.  
If $\om$ is a closed form on $M$ that is symplectic in the neighborhood $U_p$ of $p$,
then one can obtain a family of closed forms $\Tom_\eps$ on $\TM$
for small $\eps>0$, by thinking of the blow up as the manifold obtained from $(M,\om)$ by cutting out a symplectic ball in $U_p$  of  radius 
$\sqrt{\eps/\pi}$ and identifying the boundary of this ball with the exceptional divisor via the Hopf map: see~\cite{MS1}.
 
 Here is a simplified version of one of the main results in \cite{Mcbl}.
 
 \begin{prop}\labell{prop:csymp} Let $(M,a)$ be a $c$-symplectic manifold.  Then there is a  homomorphism $\Tf_*: \pi_2(M) 
  \oplus \Z\to \pi_1\bigl(\Diff\TM\bigr)$ whose kernel is 
contained in the kernel of the rational Hurewicz homomorphism
$\pi_2(M)\to H_2(M;\Q)$. Moreover, if $(M,\om)$ is symplectic, we may construct $\Tf_*$ so that  it takes $\pi_2(M)\oplus \{0\})$ into $\pi_1(\Ham(\TM,\Tom_\eps))$ for sufficiently small $\eps>0$.
\end{prop}

This result has several consequences that  throw light on 
Questions~\ref{q2} and \ref{q3}. The first 
was observed by K\c edra~\cite{K1} (and was one of the motivations for
\cite{Mcbl}.).

\begin{cor}[\cite{K1}]  There are symplectic $4$-manifolds $(M,\om)$ such that $\pi_1(\Ham\, M)$ is nonzero, but that do not support any $S^1$ action.
\end{cor}
\begin{proof}  There are many blow up manifolds  that do not admit circle actions.   For example, in dimension $4$ Baldridge~\cite{Bald}
 has shown that if $X$ has $b^+>1$ and admits a circle action with a fixed point then its Seiberg--Witten invariants must vanish.  Since 
 manifolds that admit fixed point free circle actions must have zero Euler characteristic, this implies that no simply connected K\"ahler surface with $b^+>1$  admits a circle action.  Thus the blow up of a $K3$ surface 
 has no circle actions, but, by Proposition \ref{prop:csymp}, does have nontrivial
 $\pi_1(\Ham)$ and 
$\pi_1(\Diff)$.  \end{proof}

When $M$ is symplectic rather than $c$-symplectic, 
the map $\Tf_*$  of Proposition~\ref{prop:symp} takes the factor 
$\pi_2(M)$ into $\pi_1(\Ham)$.  However, the elements 
in the image of the $\Z$ factor 
are constructed using a twisted blow up construction 
 and so most probably do not
 lie in $\pi_1(\Ham)$ for any $\om$.  
 When $\dim M = 4$
there are some cases when one can actually prove this.

\begin{cor}\labell{cor:diff}   Let $X$ be the blow up of $\C P^2$ 
at one point or the  blow up of $T^4$ at $k$ points for some $k\ge 1$.
  Then $\pi_1(\Diff X)$ is not generated by the images of 
  $\pi_1(\Symp(X,\om))$, as $\om$ varies over
   the space of all symplectic forms on $X$.
   \end{cor}

This result should be compared with Gromov's observation in
\cite{G} that the map
$\pi_1(\Ham(M,\om))\to \pi_1(\Diff)$  is not surjective when $M = S^2\times S^2$ and  each sphere factor has the same $\om$-area.   Seidel~\cite{Sei2} extended this to 
$\C P^m\times \C P^n$.
However,  the  elements in $\pi_1(\Diff)$ that they consider are in the image of $\pi_1(\Ham(M,\om'))$ for some other  $\om'$, and so they did not  establish the stronger statement given above.

We now sketch the proof of Proposition~\ref{prop:csymp}
 in the case when $H^1(M;\R) = 0$. 
 The first task is to define the map 
 $\Tf_*:\pi_2(M)\oplus \Z\to \pi_1(\Diff\TM)$. We shall identify 
 $ \pi_1(\Diff\TM)$ with $ \pi_2(B\Diff\TM)$, i.e. for each pair $(\al,\ell)\in \pi_2(M)\oplus \Z$ we shall construct a smooth bundle 
 $P_{(\al,\ell)}\to S^2$ with fiber $\TM$, defining $\Tf_*(\al,\ell)$  to
 be the  homotopy class of the corresponding clutching function 
 $S^1\to \Diff\TM$.    
 
 First suppose that $(M,\om)$ is symplectic.
 Then the diagonal is a symplectic submanifold
 of $(M\times M,\om \oplus \om)$  and so we can construct a
 symplectic bundle
 $(\TQ,\TOm_\eps)\to M$ with fiber $(\TM, \Tom_\eps)$ by blowing up normal to the diagonal with weight $\eps$.  The bundle $\TP_{\al}\to S^2$ is then defined to be the pullback of $\TQ\to M$ by $\al:S^2\to M$:
 $$
 \begin{array}{ccc}  \TP_{\al}&\longrightarrow& \TQ\\
 \downarrow&&\downarrow\\
 S^2&\stackrel{\al}\longrightarrow& M.
 \end{array}
 $$
 Thus we get a homomorphism 
 $$
 \Tf_*:\pi_2(M)\to \pi_1(\Symp (\TM,\Tom_\eps)).
 $$
 (This construction is basically due to K\c edra~\cite{K2}. The image 
 of the map actually lies in $\pi_1(\Ham)$ since $\Tom_\eps$ has 
 the closed extension $\TOm_\eps$; cf. Lalonde--McDuff
 \cite{LM2}.)
 
 Alternatively, instead of constructing the universal model, we can start with the product 
 $
 (S^2\times M, \Om)$, where $\Om = pr_1^*(\be) + pr_2^*\om
 $
  for some area form $\be$ on $S^2$, and then blow this up along the
graph 
$$
gr_\al: = \bigl\{(z,\al(z)): z\in S^2\bigr\}
$$
 of $\al$.  To make this work we just need to choose $\be$ so that $\Om$ is symplectic near $gr_\al$ and restricts on it to a nondegenerate form.  In particular, 
 we do {\it not} need $\om$ to be symplectic, and so can take it to be
 any representative of the class $a\in H^2(M)$. 
 But then there is another degree of freedom. If $\Om$  is not symplectic everywhere, there is no longer a canonical choice of complex structure on the normal bundle $\nu_\al$ to $gr_\al$.
 In fact, there is  $\Z$s worth of choices of complex structure $J$
 on  $\nu_\al$.   
 
 We define $\Tf_*(\al,\ell)$ to be the element of $\pi_1(\Diff\TM)$ obtained by blowing up
 along $gr_\al$ with complex structure chosen 
 so that $c_1(\nu_\al,J) = \ell$.  (This only depends on $\ell$ not on
  the particular choice of $J$.)
For example, we get a \lq\lq nonsymplectic" element in $\pi_1(\TM)$ by
blowing up $S^2\times M$ along the trivial section $S^2\times \{pt\}$  using a nontrivial complex structure.  (When $M=\C P^2$ we prove
 Corollary~\ref{cor:diff} by showing that the corresponding loop is not homotopic to any loop of  symplectomorphisms.)

\begin{proof}[Proof of Proposition~\ref{prop:csymp}]
 If $\al$ is not in the kernel of the Hurewicz map, there is some class in $H^2(M)$ that does not vanish on it.  Hence, by slightly perturbing $a$ if necessary we may suppose that $\la: = a(\al)\ne 0$.
  The bundle $\TP_{\al,\ell}\to S^2$ constructed above is equipped with a closed form $\Tom_\eps$ that is symplectic near the exceptional divisor.
We show below that
\begin{equation}\labell{eq:vol}
\vol(\TP_{\al,\ell},{\Tilde\Om_\eps}) = \mu_0V -  
v_\eps\left(\mu_0 + \la -\frac\ell{n+1} \eps\right),
\end{equation}
where  $\la$ is as above, $\mu_0=\int_{S^2}\be$ is the area of the base
with respect to the chosen area form $\be$,  $V
= \frac1{n!}\int_Ma^n$ is the volume of $(M,a)$,\and  $v_\eps: 
= \frac{\eps^n}{n!}$ is the volume of a ball of capacity $\eps$. Thus 
$\mu_0V$ is the volume of $(S^2\times M,\Om)$.

Now observe  that the underlying smooth bundle $\TP\to S^2$ does not depend on the choice of $\eps$.  Therefore, if $\TP_{\al,\ell}\to S^2$
 were a smoothly trivial fibration, then, for all $\eps$, the volume of $\TP_{\al,\ell}$ 
would be the product of $V-v_\eps: = \vol(\TM,\Tom_\eps)$ with the \lq\lq size"
 $\mu$ of the base. Since $\mu$ could be measured by  integrating  $\Tilde\Om_\eps$ over a  section of $\TP$ which is the same for all $\eps$, $\mu= \mu_1 + k\eps$ would be a linear function\footnote
 {
 We cannot assume that it is independent of $\eps$ since the diffeomorphism 
 may not converge to a product as $\eps\to 0$.  For example, consider the 
(trivial) bundle obtained by blowing up  $T^2\times T^2$ along the diagonal.}
  of $\eps$.  Therefore,  the two polynomial functions 
$\vol (\TP_{\al,\ell},{\Tilde\Om_\eps})$ and $(V-v_\eps)(\mu_1 + k\eps)$ 
would have to be equal.
This is possible only if $\la =0$ and also $k=\ell=0$.  
The result follows.

  It remains to derive the formula for 
  $\vol(\TP_{\al,\ell},{\Tilde\Om_\eps})$.
   First assume that $\ell = 0$.  Then the section 
   $gr_\al$ has trivial (complex) normal bundle in $S^2\times M$. 
   Hence it has a  neighborhood $U_\eps\subset U$ that is symplectomorphic to a product $ gr_\al\times B^{2n}(\eps)$.
    Thus, 
the volume of $U_\eps$ with respect to the form $\Om$ is $\vol(U_\eps) = v_\eps(\mu_0+\la)$ where $\mu_0: = \int _{S^2}\be$.  (Recall that $\be$
is an arbitrary area form on the base $S^2$.)
Since we construct  the blow up by cutting out
$U_\eps$ from  $(P, \Om)$ and identifying the boundary via the Hopf map,  we have 
$$
\vol(\TP_{\al,0},{\Tilde\Om_\eps} ) = \mu_0V -  v_\eps(\mu_0 + \la),
$$
as claimed.

Now consider the case when $\ell \ne 0$.  
Then  the normal bundle to
$gr_\al$ in $M\times S^2$ is isomorphic to the product $\C^{n-1}\oplus L_{\ell}$, where $L_\ell\to S^2$ is the holomorphic line bundle with $c_1=\ell$.  Therefore, we can choose $\Om$ so that it restricts in some neighborhood of $gr_\al$ to the product of 
a ball in $\C^{n-1}$ with 
 a $\de$-neighborhood  $\Nn_\de(L_\ell)$ 
of the zero section of $L_\ell$.  Identifying $\Nn_\de(L_\ell)$ with
part of the $4$-dimensional symplectic toric manifold $\PP(L_\ell\oplus\C)$, we can see that its volume is 
$h\de - \ell \de^2/2$, where $h = $ area of zero section and
$\de = \pi r^2$ is the capacity (or area) of the disc normal to $gr_\al$.
(Recall that this volume is just the area of a small neighborhood of the appropriate edge  of the moment polygon.)
Therefore, since $h = \mu_0+\la$ here, 
\begin{eqnarray*}
\vol(U_\eps) &=& \int_0^{\sqrt{\frac \eps\pi}} \vol(S^{2n-3}(r))\cdot\vol 
(\Nn_{\eps-\pi r^2}(L_\ell))\,dr\\
& = & \frac{(\mu_0+\la)\eps^{n}}{n!} - \frac{\ell\eps^{n+1}}{(n+1)!}\\
& = & v_\eps\bigl(\mu_0+\la-\frac{\ell}{n+1} \eps\bigr).
\end{eqnarray*}
Everything in the previous calculation remains valid except that we have to add $ \frac{\ell}{n+1}\eps v_\eps $ to the volume of $\TP_{\al,\ell}$.  This completes the proof.
 \end{proof}

\section{Quantum methods}\labell{s:quant}

By quantum methods we really mean the use of $J$-holomorphic curves and quantum homology.  
The main tool from this theory that has been used to understand $\pi_1(\Ham)$ is the Seidel representation.  Readers may consult the survey articles ~\cite[\S2.4]{Mcgr} or \cite{Mcint}
for a brief introduction to $J$-holomorphic curves,
or the book \cite{MS1} for a more detailed presentation.

\subsection{The Seidel representation}\labell{ss:seid}

Consider the small quantum homology  $QH_*(M):= H_*(M)\otimes \La$ of $M$.  
 Here $\La: = \La^{univ}[q,q^{-1}]$ where $q$ is a polynomial
 variable of degree $2$
 and  $\La^{univ}$ denotes the
generalized Laurent series ring with elements 
$\sum_{i\ge 1} r_it^{\ka_i}$, 
where $r_i\in \Q$ and $\ka_i\in \R$ is a strictly
 decreasing sequence that 
tends to $-\infty$.
 We write the elements of
 $QH_*(M)$ as infinite sums $\sum_{i\ge 1} a_i\otimes q^{d_i}t^{\ka_i}$, where $a_i\in H_*(M;\Q)=:H_*(M)$,  $|d_i|$ is bounded  and $\ka_i$ is as before.
 The term $a\otimes q^{d}t^{\ka}$ has degree $2d+\deg a$.
 
 The quantum product $a*b$ of the elements $a,b\in 
 H_*(M)\subset QH_*(M)$ is defined as follows.  Let $\xi_i, i\in I$, be a basis for $H_*(M)$ and write $\xi_i^*, i\in I,$ for the basis of $H_*(M)$ that is dual with respect to the 
 intersection pairing,  that is $\xi_j^*\cdot\xi_i = \de_{ij}$.  Then 
\begin{equation}\labell{eq:QH}
 a*b: = \sum_{i,\be \in H_2(M;\Z)} \bla a,b,\xi_i\bra^M_\be\,\, \xi_i^*\,\otimes q^{-c_1(\be)} t^{-\om(\be)},
\end{equation}
 where $\bla a,b,\xi_i\bra^M_\be$ denotes the Gromov--Witten invariant in $M$ that counts curves in class $\be$ through the homological
 constraints
 $a,b,\xi_i$.
 Note that if $(a*b)_\be: = \sum_i \bla a,b,\xi_i\bra^M_\be \xi_i^*$, then $(a*b)_\be\cdot c = \bla a,b,c\bra^M_\be$. Further, 
 $\deg(a*b) = \deg a + \deg b - 2n$, and the identity element is $\1: = [M]$.
 The product is extended to $H_*(M)\otimes \La$ by linearity over $\La$.

The {\it Seidel representation} is a homomorphism $\Ss$ from $\pi_1(\Ham(M,\om))$ to the degree $2n$ multiplicative units $QH_{2n}(M)^\times$ 
of the small quantum homology ring first considered by 
Seidel in~\cite{Sei}.  One way of thinking of it (which corresponds to the formula given below) is to say that it  "counts" all the sections of the bundle $P_\ga\to S^2$ associated to $\ga$.   However, it can also be considered
as the Floer continuation map around the loop $\ga$ in the
Hamiltonian group, which is Seidel's original point of view.  This second point of view makes it more clear why it is a homomorphism; the connection between them is discussed in Lalonde--McDuff--Polterovich~\cite{LMP}. 

To define $\Ss$, observe that  each loop $\ga=\{\phi_t\}$ in
  $\Ham$ gives rise to an 
$M$-bundle
$P_\ga\to S^2$ defined\footnote
{Different papers
 have different sign conventions.  Here we use those of
\cite[\S2.1]{MT}.  Thus we define the Hamiltonian vector field $X_H$  of
$H_t$ by the equation $ \om(X_t,\cdot) = -dH_t$.  If $K:S^2\to \R$ is the height function and we take the form $dx_3\wedge d\theta$ on $S^2$ then
the induced $S^1$ action is generated by $\p/\p\theta$.  Moreover the signs  have been chosen so that
$K$ has positive weights at its minimum and negative weights at its maximum.}
 by the clutching function $\ga$:
$$
P_\ga: =  (D_+\times M)\;\cup\;  (D_-\times M)/\!\sim, \;\;\;\;(e^{2\pi it},\,\phi_t(x))_+\sim
(e^{2\pi it},\, x)_-.
$$
where $D_{\pm}$ are copies of the unit disc in $\C$.
Because the loop $\ga$ is Hamiltonian, the fiberwise symplectic form $\om$ extends to a closed form $\Om$ on $P_\ga$, that we can arrange to be symplectic by adding to it the pullback of a suitable form on the base $S^2$.

In the case of a circle action with normalized moment map $K:M\to \R$ we may simply take
$(P_\ga, \Om)$ to be the quotient $(S^3\times_{S^1} M, \Om_c)$, where $S^1$ acts diagonally on $S^3$ and on its product with $M$ and $\Om_c$ pulls back to
$\om + d\bigl((c-K)\al\bigr)$.  Here $\al$ is the standard contact form on $S^3$ normalized so that it descends to an area form on $S^2$ with total area $1$, and $c$ is any constant larger than the maximum $K_{\max}$ of $K$.
(This last condition implies that $\Om_c$ is nondegenerate.) 
Points $p_{\max}, p_{\min}$  in the fixed point sets $F_{\max}$ and $F_{\min}$ give rise to sections
$s_{\max}: = S^2\times \{p_{\max}\}$ and $s_{\min}:=S^2\times  \{
p_{\min}\}$.  Note that 
our orientation conventions are chosen so that the integral 
of $\Om$ over the section $s_{\min}$ is {\it larger} than that over $s_{\max}$. For example, if $M=S^2$ and $\ga$ is a full rotation,
$P_\ga$ can be identified with the one point blow up of $\C P^2$, and $s_{\max}$ is the exceptional divisor. In the following we denote particular sections as $s_{\max}$ or $s_{\min}$, while writing $\si_{\max},\si_{\min}$ for the homology classes they represent.

The bundle $P_\ga\to S^2$ carries two canonical cohomology classes, the first Chern class $c_1^{\Ver}$ of the vertical tangent bundle and
the coupling class $u_\ga$, the unique class that extends the
fiberwise
symplectic class $[\om]$ and is such that $u_\ga^{n+1}=0$.
Then we define
\begin{equation}\labell{eq:S}
\Ss(\ga): = \sum_{\si, i} \bla \xi_i\bra^P_\si\,\xi_i^*\otimes q^{-c_1^{\Ver}(\si)} 
t^{-u_{\ga}(\si)}\in H_*(M)\otimes \La,
\end{equation}
Cf. \cite[Def.~11.4.1]{MS2} and \cite{Mcq}. Here $\bla \xi_i\bra^P_\si$
denotes the Gromov--Witten invariant that
counts curves in class $\si$ though the single constraint $\xi_i$ that we imagine represented in a fiber of $P$; similarly,
below $\bla b,\xi_i\bra^P_\si$ counts curves through two fiberwise constraints.  The sum is taken over all section classes in $P$ and over 
the basis $\xi_i$ for $H_*(M)$. 

Further for all $b,c\in H_*(M)$
\begin{equation}\labell{eq:Sa1}
\Ss(\ga)*b = \sum_{\si,i} \bla b,\xi_i\bra^P_\si\,\xi_i^*\otimes q^{-c_1^{\Ver}(\si)} 
t^{-u_{\ga}(\si)}.
\end{equation}

In  general it is very hard  to calculate $\Ss(\ga)$.  
The following result is essentially due to Seidel~\cite{Sei}; see also McDuff--Tolman~\cite{MT}.   We denote by $F_{\max}$ (resp. $F_{\min}$) the fixed point component
on which the normalized Hamiltonian $K:M\to \R$ that generates the action takes its maximum (resp. minimum), and by $K_{\max}$
(resp. $K_{\min}$) the maximum (resp. minimum) value of $K$.

\begin{prop}\labell{prop:max}  Suppose that  the Hamiltonian
 circle action $\ga$ is generated by the normalized Hamiltonian $K$. Then 
$$
\Ss(\ga): = 
a_{\max}\otimes q^{m_{\max}} 
t^{K_{\max}} + \sum_{\be\in H_2(M;\Z),\;\om(\be)>0} a_\be\otimes q^{m_{\max}-c_1(\be)} 
t^{K_{\max}-\om(\be)},
$$
where $m_{\max}: = -c_1^{\Ver}(\si_{\max})$ and
$a_{\max}\in H_*(F_{\max})\subset H_*(M)$.  Moreover, if the action is locally semifree near $F_{\max}$ then $a_{\max} = [F_{\max}]$.
\end{prop}

The first statement below is an immediate consequence of
the above result, but the second requires some work.
 
\begin{cor}[\cite{MT}]\labell{cor:seid}  Let $\ga$ be a 
Hamiltonian circle action on $M$.\MS

\NI
{\rm (i)}  If the  action  is semifree near either $F_{\max}$ or $F_{\min}$ then $\ga$ represents a nonzero element in $\pi_1(\Ham)$.\MS

\NI {\rm (ii)}  Suppose that $\ga$ is contractible in $\Ham$, and choose $k$ so that for all $p\in M$ the order of the stabilizer subgroup at $p$  is at most $k$.
Then, after replacing $\ga$ by $\ga^{-1}$ if necessary,
$$
K_{\max}\le |K_{\min}|\le (k-1)K_{\max}.
$$
In particular, if no stabilizer subgroup has order $> 2$ then $K_{\max} = 
-K_{\min}$.
\end{cor}

Before leaving this topic, we remark that the Seidel representation can
sometimes  be used to give information about $QH_*(M)$ itself.  
Tolman--Weitsman~\cite{TW} show that if $M$ supports a semifree $S^1$ action with isolated fixed points then there is a natural isomorphism between the equivariant cohomology of $M$ and that of a product of $2$-spheres.
In \cite{Gonzalez,Gonzalez2} Gonzalez  used the Seidel representation to improve this, showing that:

\begin{quote} {\it Under these hypotheses $M$ has the same  quantum homology as a product of $2$-spheres. Moreover, if $\dim M \le 6$ then $M$ is equivariantly symplectomorphic to a product of $2$-spheres.}
\end{quote}

\subsection{Applications to Hofer geometry}\labell{ss:Hof}

We give a very brief introduction to Hofer geometry.  For more details see \cite{H,Pbk,MS2,Mcv}.
The {\it Hofer length} $\Ll(\{\phi^H_t\})$ of a Hamiltonian path $\phi^H_t$ 
from $\phi^H_0=id$ to $\phi^H_1: = \phi$ is defined to be
$$
\Ll(\{\phi_t\}) = \int_0^1 \Bigl(
\max_{x\in M}H_t(x) - \min_{x\in M}H_t(x)\Bigr)dt.
$$
The Hofer (pseudo)norm of an element $\Tphi = (\phi, \{\phi_t\})$ in 
the universal cover $\THam$
of $\Ham$ is then defined to be the infimum of the lengths of the paths 
to $\phi$ that are homotopic to $\{\phi_t\}$.  

\begin{quest}  Is the Hofer pseudonorm a norm?
\end{quest}

It is known that this pseudonorm descends to a norm $\|\cdot\|$ on $\Ham$, i.e. $\|\phi\|=0\Leftrightarrow \phi=id$.   (It is obvious that
$\|id\|=0$ but the other implication is quite hard, see \cite{H,LM}.) Therefore this question is equivalent to the following:

\begin{quest}  Are there any nontrivial elements $\ga\in  \pi_1(\Ham)$ whose Hofer length $\|\ga\|$ is zero?
\end{quest}

Since the Hofer norm is conjugation invariant, the corresponding metric
given by
$$
\rho(\phi,\psi): = \|\phi\psi^{-1}\|
$$
is biinvariant.  There are many interesting open questions in Hofer geometry.  Here we focus on the question of whether the diameter of $\Ham$ can ever be finite.

One way to estimate the length of a Hamiltonian path is to use the Schwarz--Oh {\it spectral invariants}; see \cite{Sch,Oh,U} and \cite[Ch~12.4]{MS2}.
For each element $\Tphi\in \THam$ and each element $a\in QH_*(M)$
one gets a number $c(a,\Tphi)\in \R$ with the following properties:
\begin{eqnarray}\labell{eq:1}
&&-\|\Tphi\| \le c(a,\Tphi)\le \|\Tphi\|;\\\labell{eq:2}
&& c(\Tphi, \la a) =  c(a,\Tphi) +\nu(\la)\;\;\mbox{ for all }\la\in \La,\\\labell{eq:3}
&& 
c(a,\Tphi\circ\ga)= c(\Ss(\ga)*a,\Tphi)\;\;\mbox{ for all } \ga\in 
\pi_1(\Ham),\\\labell{eq:4}
&& 
c(a*b,\Tphi\Tpsi)\le  c(a,\Tphi) + c(b,\Tpsi)\;\;\mbox{ for all } \Tphi,\Tpsi \in \THam,\;a,b\in QH_*(M).\end{eqnarray}
In ({\ref{eq:2}), $\nu:QH_*(M)\to\R$ is the valuation given by
$$
\nu(b): = \max \{-\ka_i|b_i\ne 0\},\quad\mbox{where } 
b = \sum b_iq^{d_i}t^{-\ka_i},\;b_i\in H_*(M).
$$
These invariants are defined by looking at the filtered Floer complex\footnote
{
The differential in this complex  depends on the choice of a suitable
family of almost complex structures $J_t$, but the spectral values are  independent both of this choice and of the choice of $H$.} 
$CF_*(H,J)$ of the generating 
(normalized) Hamiltonian $H$, and by \cite{Oh1,U} turn out to be particular critical values of the corresponding
action functional $\Aa_H$.\footnote
{
This is essentially the same as the action functional used in \S\ref{ss:act}.  However, because  we want it to take values in $\R$ we must now define it on
the universal cover $\Tilde{\Ll_0M}$ of the loop space $\Ll_0M$.}  Thus each invariant $c(a,\Tphi)$ corresponds to a particular fixed point of the endpoint $\phi_1\in \Ham$ of $\Tphi$.  Property ({\ref{eq:3}) above explains how they depend on the path $\Tphi$. 

It is usually very hard to calculate these numbers. 
However, if
 $\Tphi^H$ is generated by a $C^2$-small  mean normalized Morse function $H:M\to \R$, then the invariants $c(a,\Tphi^H)$ for $a\in H_*(M)$
 are the same as the corresponding invariants $ c_M(a,\Tphi^H)$  obtained from the Morse complex $CM_*(H)$ of $H$.  These are defined as follows.
For each $\ka\in \R$, denote by $CM_*^\ka(M,H)$ the subcomplex of the Morse complex generated by the critical points $p$ with critical values $H(p)\le \ka$.
 Denote by $\io_\ka$ the inclusion of the homology  $H_*^\ka$
 of this subcomplex  into $H_*^{\infty} \cong H_*(M)$. Then 
 for each $a\in H_*(M)$ 
\begin{equation}\labell{eq:M}
 c_M(a,\Tphi^H): = \inf\{\ka: a\in {\rm Im}\,\io_\ka\}.
\end{equation}
Moreover, in this case the filtered Floer complex $CF_*(H,J)$ is simply the tensor product $CM_*(H)\otimes \La$ with the obvious product filtration.
Hence
\begin{equation}\labell{eq:M1}
 c\bigl(\sum a_i\otimes q^{d_i} t^{\ka_i},\Tphi^H\bigr) = \sup_i \,(c_M(a_i,\Tphi^H) + \ka_i) 
 \end{equation}

As pointed out by Ostrover~\cite{Os}, one can use the continuity properties of the $c(a,\Tphi)$ to prove that $\THam$ has infinite diameter. Indeed,
let $H$ be a small Morse function as above, choose an open set $U$ 
that is displaced by $\phi_1^H$ (i.e. $\phi_1^H (U)\cap U=\emptyset$)
and let $F:M\to U$ be a function with support in $U$ and with nonzero integral $I: = \int_M F\om^n$.  Denote the flow of $F$ by $f_t$ and consider the path $\Tphi_s: = \{f_{ts}\phi_t^H\}_{t\in [0,1]}$ for $s\to \infty$.  This is generated by the Hamiltonian
$$
F_s\# H: = F_s + H\circ f_{st}.
$$
The corresponding normalized Hamiltonian
is
$$
K_s: = F_s + H\circ f_{st} - sI.
$$ 

By construction, $f_{s}\phi_1^H$ has the same fixed points
as $\phi_1^H$, namely the critical points of $H$.  Hence the continuity and spectrality  properties of $c(a,\Tphi_s)$ imply that for each $a\in QH_*(M)$ the fixed point $p_a$ whose critical value is
$c(a,\Tphi_s)$  remains unchanged as $s$ increases.
But the spectral value does change. In fact, if $a\in H_*(M)$, then when $s=0$  there is a critical point $p_a$ of $H$ such that
$c(a,\Tphi_0) = c_M(a,\phi_1^H) = H(p_a).$  Hence 
\begin{equation}\label{eq:sp}
c(a,\Tphi_s) = K_s(p_a) =  H(p_a) - sI,\;\;\mbox{ for all } s\in \R,\; a\in H_*(M).
\end{equation}
By ({\ref{eq:1}) it follows that $\THam $ has infinite Hofer diameter. 

Schwarz proved the following result in ~\cite[\S4.3]{Sch}
 under the assumption that  both $\om$ and $c_1$ vanish on 
$\pi_2(M)$.  However, his argument (which  
is based on an idea due to Seidel)
works equally well if one just assumes that $\om$ vanishes.

\begin{prop} If $\om$  vanishes on 
$\pi_2(M)$ then every Seidel element  has the form
$$
\Ss(\ga) = (\1+x)\otimes \la,\quad\mbox{where }\; \nu(\la) = 0,\;x\in H_{<2n}(M)\otimes\La.
$$
Moreover the spectral invariants descend to $\Ham$.  
In particular, $\Ham$ has infinite diameter. 
\end{prop}
\NI {\it Sketch of proof.}
If $\om=0$ on $\pi_2(M)$, quantum multiplication is the same as the intersection product. Hence all units in $QH_*(M)$ have the form
 $\1\otimes \la + x$ where $\la\in\La$ is nonzero and $x\in H_{<2n}(M)
 \otimes\La$.  In particular there must be a section class $\si_0$  of $P_\ga$
  with $c_1^{\Ver}(\si_0)=0$ that contributes to the coefficient of $\1$.
Moreover, because all sections of $P_\ga$  have the same energy,  
the Seidel element has the form $\Ss(\ga) = (\1 + x)\otimes t^\ka$
where $x\in H_{<2n}(M)[q,q^{-1}]$. Therefore it suffices to show that 
$\ka: = -\nu(\si_0) = 0$.

Now choose a generic almost complex structure $J$ on $P$ that is compatible with the fibration $P\to S^2$ and consider the moduli space $\Mm$ of $J$ holomorphic sections in class $\si_0$, parametrized as sections.   This space is always compact because there can be no fiberwise bubbles, and is a manifold for generic $J$.  Hence there is a commutative diagram
$$
\begin{array}{ccc} S^2\times \Mm&\stackrel{ev}\longrightarrow& P\\
\downarrow&&\downarrow\\
S^2&\stackrel{=}\longrightarrow& S^2,\end{array}
$$
where the top horizontal map is the evaluation map.  Since  $ev$ maps each fiber $\{z\}\times \Mm$  to the corresponding fiber $\pi^{-1}(z)$ with   positive degree, the identity $(u_\ga)^{n+1}=0$ implies that the class
$\ev^*(u_\ga)\in H^2(S^2\times \Mm)$ is pulled back from $\Mm$. Hence  
$$
\ka : =- \int_{\si_0} u_\ga  \;=\; \int_{S^2\times \{pt\}}\,ev^*(u_\ga) \;= 0.
$$

This proves the first statement.   It follows  that
 $\nu(\Ss(\ga))\le 0$ for all $\ga$.  But it is impossible to have
 $\nu(\Ss(\ga))<0$ for some $\ga$ 
 because $\Ss(\ga)*\Ss(\ga^{-1})=\1$ and  
$\nu(a*b)\le \max(\nu(a),\nu(b))$.    Hence $\nu(\Ss(\ga))=0$ for all $\ga$.
But then the multiplicative relation ({\ref{eq:4}) implies that
$$
c(a,\Tphi\circ\ga)\le c(a,\Tphi) + c(\1,\ga) = c(a,\Tphi) + \nu(\Ss(\ga)) =c(a,\Tphi), 
$$
for all $\ga\in \pi_1(\Ham).$ But then $c(a,\Tphi)=c(a,\Tphi\circ\ga\circ\ga^{-1})\le c(a,\Tphi\circ\ga)$.  Hence $c(a,\Tphi)=c(a,\Tphi\circ\ga)$ for all $\ga$, and so
 the spectral invariants descend.
\QED
%
%




These arguments are significantly extended in \cite{Mcmdr}.  For example, 
if $\dim H_2(M)>1$, the condition $[\om]=0$ on $\pi_2(M)$ in the above proposition can be replaced by the much weaker condition that all the nontrivial genus zero $3$-point Gromov--Witten invariants of $M$ vanish.

 Here is a somewhat different approach to this question that gives information
 about the Hofer diameter even when the spectral invariants do not descend.  
We shall say that $\pi_1(\Ham)$ is  {\it spectrally asymmetric} if
$$
\liminf_{\|\ga\|\to \infty} \frac{|\nu(\Ss(\ga))+\nu(\Ss(\ga^{-1}))|}{\|\ga\|} \ge \eps> 0.
$$
Roughly speaking, this means that for  loops $\ga$ that are long with respect to the Hofer norm $\|\ga\|$, their spectral \lq\lq norms"
$\nu(\Ss(\ga))$ and $\nu(\Ss(\ga^{-1}))$ are not comparable.
For example, if $\pi_1(\Ham)\cong\Z$ with generator $\al$, this 
holds if $$
 \lim_{k\to\infty} \frac{\nu(\Ss(\al^k))}k \ne  -
 \lim_{k\to\infty} \frac{\nu(\Ss(\al^{-k}))}k 
 $$
because $\|\al^k\|\le k\|\al\|$.

\begin{lemma}\labell{le:inf} $\Ham$ has infinite diameter
whenever $\pi_1(\Ham)$ is   spectrally asymmetric. 
\end{lemma}
\NI {\it Sketch of proof.}  
Consider the path $\Tphi_s$ of equation (\ref{eq:sp}).  Then by equation
(\ref{eq:M1}) there is
for any $\ga\in \pi_1(\Ham)$  a critical point $p_\ga$ of $H$ such that 
\begin{eqnarray*}
c\bigl(\1,\ga \Tphi_s\bigr)
&=& c\bigl(\Ss(\ga),\Tphi_s\bigr) \\&=& c\bigl(\Ss(\ga),f_s\phi^H_1\bigr)\\
& = &
c\bigl(\Ss(\ga),\phi^H_1\bigr) - sI \\
&=& H(p_\ga)  + \nu(\Ss(\ga)) - sI.
\end{eqnarray*}
Now suppose that $\Ham$ has diameter $D$.
 Because the Hofer norm is symmetric, i.e. $\|\phi\| = \|\phi^{-1}\|$, for each $s$ there must be a loop $\ga_s$ such that 
$$
|c\bigl(\1,\ga_s \Tphi_s\bigr)|\le D,\quad |c\bigl((\1,(\ga_s\Tphi_s)^{-1}\bigr)|\le D.
$$
Therefore, 
$$
|c\bigl(\Ss(\ga_s),\phi^H_1\bigr) - sI|\le D, \quad
|c\bigl(\Ss(\ga_s^{-1}),\phi^{-H}_1)\bigr) + sI|\le D.
$$
But by (\ref{eq:1})  $|\nu(\Ss(\ga_s))| = |c(\1,\ga_s)|\le \|\ga_s\|$.
Hence  $\|\ga_s\|\to \infty$ as $s\to \infty$, and so  
\begin{equation}\labell{eq:as}
\lim_{s\to \infty}\frac 
{|c\bigl(\Ss(\ga_s),\phi^H_1\bigr) + c\bigl(\Ss(\ga_s^{-1}),\phi^{-H}_1\bigr)|}
{\|\ga_s\|} =0.
\end{equation}
But we saw above that there are critical points $p,q$ of $H$ such that
\begin{eqnarray*}
&& c\bigl(\Ss(\ga_s),\phi^H_1\bigr)=H(p) + \nu (\Ss(\ga_s)),\\
&& c\bigl(\Ss(\ga_s^{-1}),\phi^{-H}_1\bigr)= -H(q) + \nu (\Ss(\ga_s^{-1})).
\end{eqnarray*}
 Hence the LHS of  (\ref{eq:as}) has the same limit as
 $$
 \frac{|\nu (\Ss(\ga_s)) +\nu (\Ss(\ga_s^{-1}))|}{\|\ga_s\|},
 $$
 and this is $\ne 0$ by assumption.  This contradiction shows that 
 $\Ham$ cannot have finite diameter.
\QED

Now consider the one point blow up $M$ of $\C P^2$ with symplectic form $\om_a$ that takes the value $\pi a$ on the exceptional divisor and $\pi$ on the line. 
$(M,\om_a)$ is a \lq\lq small" blow up if $3a^2<1$, i.e. the exceptional divisor is smaller than it is in the monotone case.

\begin{cor}  $\Ham(M,\om_a)$ has infinite diameter when $3a^2<1$.
\end{cor}
\begin{proof}  By \cite{AM} $\pi_1(\Ham(M,\om_a))= \Z$.   
By  \cite[Prop.~5.3]{Mcv}, it has a
 generator  $\al$ with
$\Ss(\al) = Q\otimes t^{\mu}$ where $\mu =
\pi\frac{(1-a^2)^2}{12(1+a^2)}$.  Moreover when $3a^2<1$ there are positive constants $c_1, c_2$ such that
$$
\nu(Q^{-k})\ge c_1 k,\quad \nu(Q^{k})\ge c_2 k,\qquad k\ge1.
$$ 
Hence, because $\|\al^k\|\le k\|\al\| $,
$$
\liminf_{k\to \infty} \frac{|\nu(\Ss(\al^k))+
\nu(\Ss(\al^{-k}))|}{\|\al^k\|} \;\ge\;
\liminf_{k\to \infty} \frac{|(c_1 -\mu)k + (c_2+\mu)k|}{k\|\al\|} \;>\;0.
$$
Thus $\Ham(M, \om_a)$ is spectrally asymmetric.
\end{proof}

It is shown in \cite{Mcv}  that  $|\nu(Q^k)|$ is bounded when 
$3a^2\ge1$.  Therefore  in this case our methods give no information on the Hofer diameter.

 \begin{rmk}\rm  In the above Lemma~\ref{le:inf} we exploited the symmetric definition of the
  Hofer norm.  As pointed out in \cite{Mcv} there are
   variants of this norm with more 
  asymmetric definitions; for example one can measure the size of $\phi\in \Ham$ by separately minimizing the positive and negative parts of the Hofer lengths over all paths in $\Ham$ from $id$ to $\phi$.   This seminorm is not known to be nondegenerate for all $M$, though it is
   for the one point blow up of 
  $\C P^2$.   Moreover the argument in  Lemma~\ref{le:inf} would not apply to it.
  \end{rmk}

\subsection{The Calabi homomorphism}

Finally, we show how to calculate the Calabi homomorphism
$\pi_1(\Ham^c) \to \R$ for certain loops on  noncompact manifolds $M$,
and also explain its relation to the Seidel element.

Suppose that the loop $\ga\in \pi_1(\Ham^c\,M)$ is supported in an open subset $U\subset M$ such that $(U,\om)$ 
 can be symplectically embedded into some closed
manifold $(X,\om_X)$.  Consider the corresponding bundle
$X\to P_{\ga_X} \to S^2$ and the corresponding Seidel element 
$\Ss(\ga_X)\in QH_*(X)$.   Then $P_{\ga_X}$ has a family of flat sections
$S^2\times \{p\}, p\in X\smallsetminus U.$ 
Suppose now that  these sections contribute nontrivially to the coefficient
of $\1\otimes t^\ka$ in $\Ss(\ga_X)$.  

\begin{lemma}\labell{le:cal}
$
\Cal(\ga) =  -\ka\,\Vol(X,\om_X).
$
\end{lemma}
\NI {\it Sketch of Proof.} Let $\ga_X=\{\phi_t\}$ and suppose that
 $H_t: X\to \R$ is the corresponding
 mean normalized Hamiltonian.
  Then each $H_t$ equals  some constant $c_t$ outside $U$ so that 
\begin{eqnarray*}
\Cal(\ga): &=&\frac 1{n!}\int_0^1\Bigl(\int_U(H_t-c_t)\,\om^n \Bigr)\,dt\\
&=&\frac 1{n!}\int_0^1\Bigl(\int_X(H_t-c_t)\,\om_X^n \Bigr)\,dt\\
&=& - \Vol(X,\om_X) \int_0^1 c_t\,dt.
\end{eqnarray*}
On the other hand, an easy calculation shows that
the coupling class $u_{\ga_X}$ is represented by a form that we can take to be
$pr^*(\om)$ on $D_-\times X$ and
$pr^*(\om) - d\bigl(\be(r) H_t dt\bigr)$ on $D_+\times X$ for some
 cut off function $\be$.  Then $u_{\ga_X}$ restricts to zero on 
 $ D_-\times \{p\}$ and so by the definition of $\Ss(\ga)$:
\begin{eqnarray*}
-\ka: = \int_{S^2\times\{p\}} u_\ga &=
 &-\int_{D_+\times \{p\}} d(\be(r) H_t)\, dt\\
&=& -\int_{D_+} \be'(r)dr H_t(x_0)\, dt\\
&=&- \int_{S^1} H_t(x_0) \,dt\; =\; -\int_0^1 c_t\,dt.
\end{eqnarray*}
The result follows.\QED

\begin{rmk}\rm (i)
In constructing explicit examples where $\Cal(\ga)$ does not vanish, one can proceed the other way around.  Suppose given a loop $\ga_X$ for which
the corresponding bundle  $P_{\ga_X}\to S^2$ has a section $s_X$ with trivial normal bundle and with $u_{\ga_X}(s_X) = -\ka\ne0$.
(One way to find such a section is to show that 
the term $\1\otimes t^\ka$ appears with nontrivial coefficient in
$\Ss(\ga_X)$.)  
Then we may trivialize the bundle $P_{\ga_X}\to S^2$
near $s_X$ and  take $\ga$ to be the loop in $\Ham^c(X\less \{pt\})$ 
corresponding to the fibration $\bigl(P_{\ga_X}\less s_X\bigr)\to S^2$. 
If $\Ss(\ga) = \1\otimes \la + x$ where $\la = \sum r_it^{\eps_i}$ has infinitely many nonzero terms then one would get infinitely many different loops in $\Ham^c(X\less \{pt\})$.  However,  
there would  be algebraic relations between them stemming from the long exact homotopy sequence of the fibration 
$$
\Ham^{id}(X,p)\longrightarrow \Ham\, X\stackrel{ev}\longrightarrow  {\rm Fr}\,(X),
$$
where ${\rm Fr}\,(X)$ is the symplectic frame bundle of $X$ and
$\Ham^{id}(X,p)$ denotes the subgroup of $\Ham\,X$ consisting of elements $\phi$ that fix the point $p$ and have $d\phi_p=id$.  (A standard Moser type argument shows that the inclusion $\Ham^c(X\less \{p\})\to \Ham^{id}(X,p)$ is a homotopy equivalence.)
\MS

\NI (ii)
The manifolds $M=X\less \{pt\}$ obtained in (i) have finite volume and are
 concave at infinity.  However, if one started with a bundle $\TM\to \TP\to S^2$ formed by blowing up the product
  bundle along the section $gr_\al$ then
 this bundle is trivial over $\Phi^{-1}(S^2\times (M\less V))$, where $\Phi:\TP\to P=S^2\times M$ is the blowdown map and $V$ is a neighborhood of ${\rm im\, }\al$ in $M$.     Hence we could get examples in
 $\pi_1(\Ham^c\, W)$ where $W$ is any symplectic extension of $(V,\om)$.
 Note also that 
  one can of course calculate $\Cal(\ga)$
for these blow up bundles  without 
 any reference to the Seidel element.  All one needs is the formula 
  for the coupling class given in the proof of Lemma~\ref{le:inf}.
\end{rmk}

\begin{rmk}[Questions about non-compact manifolds]\rm  The above discussion of the diameter of $\Ham$ 
assumes that $(M,\om)$ is closed.  If $M$ is not closed then one can ask a similar question about $\Ham^c$.  One really should divide into two cases here, depending on whether $M$ has finite or infinite volume.\footnote
{
Observe that it is known by a rather intricate geometric argument
 that the Hofer seminorm is a norm on $\Ham^c$ for {\it every} open manifold $M$, whether nice at infinity or not; see 
\cite{LM}.}

If $\Vol\,M$ is finite, then the universal cover $\THam\,\!^c$ of 
$\Ham^c$
has infinite Hofer diameter. To see this, observe first that
the Calabi homomorphism defined earlier on $\pi_1(\Ham^c)$ extends to
$\THam\,\!^c$ by:
$$
\Cal:\;\THam\,\!^c\to \R,\quad
 \Tphi\mapsto \int_0^1\int_MH_t\om^n\,dt,
$$
where $H_t$ is any compactly supported functions that generates $\Tphi$.
Next note the obvious fact that $\Cal(\Tphi)\le \Vol(M)\,\|\Tphi\|$.

  If $\Cal$ vanishes on $\pi_1(\Ham^c)$ then it induces a homomorphism  from $\Ham^c$ to $\R$, and so 
  $\Ham^c$ must also have infinite Hofer diameter.
I do not know what happens in the general case.
 I also do not know of any methods to help  decide 
 whether the {\it kernel}  of
$\Cal$  has infinite diameter: Ostrover's argument on spectral drift does not work here because all Hamiltonians are normalized to vanish at infinity.  
Since this question is unsolved there is not much point in wondering about the corresponding subgroup of $\Ham^c$.  (Banyaga showed that this is a perfect group; hence in many ways this kernel plays the same role for noncompact $M$ as does $\Ham$ in the closed case.)  

There are other interesting open questions about  the homomorphism $\Cal:\pi_1(\Ham^c)\to \R$.   It is known to vanish when $\om$ is exact; cf. \cite[Ch~10]{MS1}.  On the other hand we saw above that it need not vanish in general.  One obvious question is whether it always has discrete image.\footnote{
Here one should no doubt restrict to $M$ of finite type, eg those that are the interior of a compact manifold $\ov M$ with boundary. But even this topological restriction on $M$ may not be enough.
Because symplectic forms exist on open manifolds under the most mild 
topological restrictions, one may need to ask that 
$\om$ extend to $\ov M$ in order to have a hope of discreteness here.
For example,  a similar question concerning the flux homomorphism was explored in \cite{Mcf}, and examples of open manifolds $M$ were constructed  for which $\Flux(\pi_1(\Symp^cM))$ is not discrete.  However
no control at infinity was imposed.  Presumably  Ono's argument in \cite{Ono} would extend if the control at infinity were strong enough.} 
 By Lemma~\ref{le:cal} 
this is closely related to the question of the possible values of $\nu(\Ss(\ga))$
for $\ga\in \pi_1(\Ham \,X)$ where $X$ is a closed manifold containing $M$.\footnote
{
One could try to define some kind of analog to the Seidel representation, though one would have to specify exactly which
 quantum homology group  is the target of this homomorphism. 
 However, it is not at all clear how such a representation 
 could help understand the questions at hand.  Lemma~\ref{le:cal} 
 suggests to me that
 in the noncompact case 
 the role of the map $\ga\mapsto\nu(\Ss(\ga))$ is played by the Calabi homomorphism.}
  (But note that the function  $\ga\to \nu(\Ss(\ga))$ need not be a homomorphism.)

Finally we remark that in the case when $M$ has infinite volume,
 there is no direct connection between the Hofer norm and the Calabi homomorphism. Therefore,  one should tackle the question of the Hofer diameter by other means, for example by seeing if some version of the energy--capacity inequality applies.  \end{rmk}

\NI{\bf Acknowledgement.}  I wish to thank A. Pelayo for his very careful reading of an earlier version of this paper, and for making many useful suggestions to improve the clarity of the exposition.  All mistakes are of course the responsibility of the author.

\end{document}